\documentclass[11pt]{article}
\usepackage{graphicx}
\usepackage{amsmath}
\usepackage{amsfonts}
\usepackage{amssymb}

\newtheorem{rem}{Remark}

\def\no{\noindent}
\def\RR{\mathbb R}

\def\CC{\mathbb C}
\def\pmatrix{ \left( \begin{array} }
\def\endpmatrix{ \end{array} \right) }

\def\aa{\alpha}
\def\bb{\beta}
\def\cc{\gamma}

\def\ttt{\theta}
\def\lam{\lambda}

\def\bfy{{\bf y}}
\def\bfc{{\bf c}}
\def\bfe{{\bf e}}
\def\bff{{\bf f}}
\def\bfeta{\mbox{\boldmath$\eta$}}
\def\bftau{\mbox{\boldmath$\tau$}}
\def\bfdel{\mbox{\boldmath$\delta$}}
\def\tro{\tilde{\rho}}
\def\troinf{\tro_\infty}
\def\rostar{\rho^*}
\def\calA{{\cal A}}

\begin{document}

\title{Blended General Linear Methods based on \\Boundary Value
Methods in the\\ Generalized BDF family}

\author{Luigi Brugnano\footnote{Dipartimento di Matematica ``U.\,Dini'',
Viale Morgagni 67/A, 50134 Firenze (Italy),\qquad\qquad
e-mail:\,luigi.brugnano@unifi.it} \and Cecilia
Magherini\footnote{Dipartimento di Matematica Applicata
``U.\,Dini'',
  Via Buonarroti 1/C , 56127 Pisa (Italy), ~
e-mail:\,cecilia.magherini@dma.unipi.it}}

\date{\small Final revised form March 28, 2009.\\
\em Dedicated to John Butcher on the occasion of his 75th
birthday}

\maketitle

\vspace{-1cm}
\begin{abstract}
Among the methods for solving ODE-IVPs, the class of General
Linear Methods (GLMs) is able to encompass most of them, ranging
from Linear Multistep Formulae (LMF) to RK formulae. Moreover, it
is possible to obtain methods able to overcome typical drawbacks
of the previous classes of methods. For example, order barriers
for stable LMF and the problem of order reduction for RK methods.
Nevertheless, these goals are usually achieved at the price of a
higher computational cost. Consequently, many efforts have been
made in order to derive GLMs with particular features, to be
exploited for their efficient implementation.

In recent years, the derivation of GLMs from particular Boundary
Value Methods (BVMs), namely the family of Generalized BDF (GBDF),
has been proposed for the numerical solution of stiff ODE-IVPs
\cite{BrTr98}. In particular, in \cite{BrMa08-1}, this approach
has been recently developed, resulting in a new family of $L$-stable GLMs of
arbitrarily high order, whose theory is here completed and fully worked-out.
Moreover, for each one of such methods, it is possible to define a
corresponding {\em Blended} GLM which is equivalent to it
from the point of view of the stability and order properties. These
{\em blended} methods, in turn, allow the definition of efficient
nonlinear splittings for solving the generated discrete problems.

A few numerical tests, confirming the excellent potential of such
{\em blended} methods, are also reported.

\medskip
\no {\bf Keywords:} Numerical methods for ordinary differential
equations, General Linear Methods, Boundary Value Methods (BVMs),
Generalized Backward Differentiation Formulae (GBDF), Blended
Implicit Methods, blended iteration.

\medskip
\no{\bf MSC:} 65L05, 65L20, 65H10, 65Y20.

\end{abstract}

\section{Introduction}
General Linear Methods, in the form introduced by Burrage and
Butcher \cite{BuBu80}, are numerical methods for solving ODE-IVPs
which, in order to advance the integration one step, require some
information from the previous step ({\em external stages}), along
with some {\em internal stages} to be computed at the current
step. Such methods are able to describe, as extreme cases, both RK
methods and LMF (the latter, when used as Initial Value Methods
(IVMs)) \cite{Bu87}. A GLM with $r$ external stages and $s$
internal stages, when applied with stepsize $h$ for solving the
IVP
\begin{equation}\label{ivp} y' = f(t,y), \qquad y(t_0) =
y_0\in\RR^m,\end{equation}

\no is usually described as
\begin{eqnarray}\label{glm} Y^{[n]} &=& h(A\otimes I_m) F(Y^{[n]})
+(U\otimes I_m)\bfy^{[n-1]}, \\ \nonumber \bfy^{[n]} &=& h(B\otimes I_m)
F(Y^{[n]}) +(V\otimes I_m)\bfy^{[n-1]},\qquad
n=1,2,\dots,\end{eqnarray}

\no where:
\begin{itemize}

\item the matrices
$A\in\RR^{s\times s}$, $U\in\RR^{s\times r}$,
$B\in\RR^{r\times s}$, $V\in\RR^{r\times r}$ characterize the method;

\item $Y^{[n]},F(Y^{[n]})\in\RR^{sm}$ are the vector of the internal stages
and the corresponding values of the function $f$, respectively;

\item $\bfy^{[n-1]},\bfy^{[n]}\in\RR^{rm}$ are the vectors of the external
stages.

\end{itemize}
Many papers have been devoted, across the years, to the
construction and the analysis of GLMs (see, e.g., the review
article \cite{Bu06}; see also
\cite{Bu93,BuCaJa97,BuJa93,BuJa96,BuJa98,BuJa04,BuJaMi97,BuWr03}).
Clearly, the efficient implementation of GLMs depends on the
properties of the matrix $A$. In particular, we shall here
consider implicit GLMs for which the external and the internal
stages coincide, i.e., $$A=B, ~ U=V \, \in\RR^{r\times
r},\qquad\mbox{and}\qquad Y^{[n]}=\bfy^{[n]} \,\in \RR^{rm}.$$

\no Moreover, matrix $A$ is nonsingular and the order of accuracy
of each entry of $\bfy^{[n]}$ is equal to $p$ (to be specified
later), which is, therefore, the global order of the method. We
mention that sometimes GLMs with approximations having the same
order have been also called ``peer methods'' (see, e.g.,
\cite{PoWeSc06}). The approach that we shall consider is the one
introduced in \cite[Chapter 11, Section 6]{BrTr98} and is based on
Boundary Value Methods (BVMs).  In more detail, the family of BVMs
which we shall consider is that of Generalized Backward
Differentiation Formulae (GBDF) \cite{BrTr96,BrTr98}. For the
efficient solution of the generated discrete problems, we shall
consider the {\em blended} implementation of such methods. This
implementation, introduced in \cite{Br00} for block implicit
methods (see also \cite{BrMa02,BrMa04,BrMa07,BrMaMu06,BrTr01}),
naturally induces an efficient splitting procedure for the
solution of the generated discrete problems. The linear analysis
of convergence for this splitting will be done according to
\cite{BrMa08}.

With these premises, the structure of the paper is the following: in
Section~\ref{bvms} Boundary Value Methods (BVMs) are briefly sketched, along
with their block form; in Section~\ref{gbdf} the family of Generalized
Backward Differentiation Formulae (GBDF), in the BVMs class, is considered,
in order to obtain $L$-stable GLMs of arbitrarily high order; in
Section~\ref{blend} the corresponding {\em Blended} GLMs are described;
Section~\ref{impl} is devoted to the implementation details of the methods;
Section~\ref{numer} contains some numerical tests; finally,
Section~\ref{concl} contains a few concluding remarks.

\section{Boundary Value Methods (BVMs) and their block form}\label{bvms}

Boundary Value Methods (BVMs) are a relatively new class of
numerical methods for ODEs based on an unconventional use of LMF.
Even though recent developments of such methods have been also
obtained (see, e.g., \cite{MaSeTr06,MaSeTr06-1}), nevertheless,
the main reference for such methods remains the book
\cite{BrTr98}. The basic idea, on which BVMs rely, can be easily
explained by just applying a $k$-step LMF, which can be described,
by using the usual notation, as \begin{equation}\label{main}
\sum_{i=0}^k \aa_iy_{n+i} = h\sum_{i=0}^k
\bb_if_{n+i},\end{equation}

\no to the standard test equation \begin{equation}\label{test} y' = \lam y,
\qquad Re(\lam)<0.\end{equation}

\no If we assume we know the first $\nu\in\{1,\dots,k\}$ values of
the discrete solution, as well as the last $k-\nu$, say
\begin{equation}\label{addval} y_0,\dots,y_{\nu-1},
\qquad\mbox{and}\qquad y_{N-k+\nu+1},\dots,y_N,\end{equation}

\no then, under suitable conditions for the method, if
$$|z_1(q)|\le |z_2(q)| \le \dots \le |z_k(q)|$$ are the zeros of the
characteristic polynomial of the method,
\begin{equation}\label{rosi}
\rho(z)-q\sigma(z), ~~\mbox{where}~~ q = h\lam, ~~\mbox{and}~~
\rho(z) = \sum_{i=0}^k \aa_i z^i, ~\sigma(z) = \sum_{i=0}^k \bb_i
z^i,\end{equation}

\no then the discrete solution is approximately given by
\begin{equation}\label{genroot}
y_n \approx c\, \left[z_\nu(q)\right]^n, \qquad\mbox{for}\quad 0\ll
n\ll N,\end{equation}

\no with the constant $c$ independent of $N$ \cite{BrTr98}. When
$\nu=k$ we have the usual way in which LMF are used; i.e., we
approximate the continuous IVP by means of a discrete IVP. We
speak, in such a case, of an {\em Initial Value Method (IVM)}. On
the other hand, when $\nu<k$ we are approximating the continuous
IVP by means of a {\em discrete BVP}. The latter defines a {\em
Boundary Value Method (BVM)} used with {\em $(\nu,k-\nu)$-boundary
conditions}. This freedom in the choice of the number $\nu$ of the
initial conditions in (\ref{addval}), allows us to overcome the
usual Dahlquist's barriers for stable LMF (see \cite{BrTr98}, for
further details).

Nevertheless, the IVP (\ref{ivp}) only provides
the initial condition $y_0$, whereas the remaining set of $k-1$
{\em additional conditions} in (\ref{addval}) are not known. However, they
can be retrieved implicitly, by introducing a suitable set of $\nu-1$ {\em
additional initial methods},
\begin{equation}\label{init}\sum_{i=0}^k \aa_i^{(j)}y_i = h\sum_{i=0}^k
\bb_i^{(j)}f_i,\qquad j=1,\dots,\nu-1,\end{equation}

\no and $k-\nu$ {\em additional final methods},
\begin{equation}\label{fine}\sum_{i=0}^k \aa_{k-i}^{(j)}y_{N-i} =
h\sum_{i=0}^k \bb_{k-i}^{(j)}f_{N-i},\qquad j=\nu+1,\dots,k,\end{equation}

\no which are independent of the {\em main method} (\ref{main}).
Moreover, in such a way, BVMs can be also implemented in block
form ({\em block BVMs} \cite{BrTr98}), which is computationally
more appealing: as a matter of fact, block BVMs have been
successfully implemented in the computational codes GAM
\cite{IaMa} and GAMD \cite{TestSet}. In more detail, the block
method, with blocksize $r$, is used in a sequential fashion by
solving, at the $n$th step, a discrete problem in the form $${\cal
A} \otimes I_m \bfy -h{\cal B}\otimes I_m \bff = \bf0,$$

\no where \begin{equation}\label{AA}\calA = \pmatrix{ccccccc}
\aa_0^{(1)} &\aa_1^{(1)} &\dots       &\dots     &\aa_k^{(1)}\\
\vdots      &    \vdots      &        &          &\vdots\\
\aa_0^{(\nu-1)} &\aa_1^{(\nu-1)} &\dots       &\dots     &\aa_k^{(\nu-1)}\\
\aa_0    &\aa_1   &\dots       &\dots     &\aa_k\\
            &\ddots &\ddots    && &\ddots\\
            &       &\aa_0 &\aa_1      &\dots  &\dots  &\aa_k\\
&&\aa_0^{(\nu+1)} &\aa_1^{(\nu+1)}    &\dots  &\dots &\aa_k^{(\nu+1)}\\
&&\vdots &\vdots        &    &      &       \vdots\\
&&\aa_0^{(k)}     &\aa_1^{(k)}  &\dots  &\dots &\aa_k^{(k)}
\endpmatrix\in\RR^{r\times r+1},\end{equation}

\no matrix ${\cal B}$ is similarly defined by replacing the
$\aa$-s with the $\bb$-s, $$\bfy =
(~y_{(n-1)r},~y_{(n-1)r+1},~\dots,~y_{nr}~)^T, \quad \bff
=(~f_{(n-1)r},~f_{(n-1)r+1},~\dots,~f_{nr}~)^T,$$

\no and $y_{(n-1)r}, f_{(n-1)r}$ are known from the previous step
(when $n=1$, $y_0$ is the initial condition for problem
(\ref{ivp})). Block BVMs allow an easier implementation of
variable stepsize, since only the stepsize $h$ of the current
``block'' can be varied. Several families of block BVMs have then
been defined (see \cite{BrTr98}, for full details).

\section{Generalized Backward Differentiation Formulae
(GBDF) and GLMs}\label{gbdf}

Let us consider the particular class of $k$-step LMF having the
polynomial $\sigma(z)$ in its simplest form, i.e.,
$$ \sigma(z) = z^j,$$
where $j\in\{0,1,\dots,k\}$, with the coefficients of the
corresponding polynomial ~$\rho(z)$~ (see
(\ref{rosi})) uni\-quely defined by imposing that the method has the maximum
possible order $p=k$. When $j=k$, one obtains the well-known family of BDF
which, however, provides $0$-stable methods only up to $k=6$ and
$A$-stable methods up to $k=2$. Nevertheless, when
$$j=\nu\equiv \lceil(k+2)/2\rceil,$$

\no and the formula is used as a BVM with $(\nu,k-\nu)$-boundary
conditions, i.e., by fixing the values (\ref{addval}) of the
discrete solution, then it turns out that the resulting methods
are stable for all values of $k$, and have been called {\em
Generalized BDF (GBDF)} \cite{BrTr96,BrTr98,AcTr02}. By the way,
when $k=1,2$, one obtains the usual first two BDF. In view of the
implementation of GBDF as GLMs, one may assume they know the
initial values of the discrete solution. On the other hand, the
final $k-\nu$ values can be retrieved implicitly by using a
suitable set of additional final methods (\ref{fine}), consisting
of LMF having the following characteristic polynomials:
$$\sigma^{(j)}(z) = z^j, \qquad \rho^{(j)}(z) = \sum_{i=0}^k \aa_i^{(j)}z^i,
\qquad j=\nu+1,\dots,k,$$

\no with the coefficients $\{\aa_i^{(j)}\}$ uniquely defined by
imposing the order $p=k$ conditions. If we assume, for the sake of
simplicity, that a constant stepsize $h$ is used, then by
introducing the vectors
$$\bfy_{new} = \pmatrix{c}y_{n+1}\\ \vdots\\y_{n+r}\endpmatrix, \qquad
\bff_{new} =\pmatrix{c} f_{n+1}\\ \vdots\\f_{n+r}\endpmatrix, \qquad
\bfy_{old} = \pmatrix{c} y_{n-r+1}\\ \vdots\\y_n\endpmatrix,$$

\no and the $r\times r$ matrices $\calA_1$ and $\calA_2$ such that
\begin{equation}\label{A1A2}
[\calA_1\,|\,\calA_2] = \pmatrix{ccccc|cccc}
&\aa_0       &\dots       &\dots    &\dots &\dots     &\aa_k\\
&            &\ddots &      &       &  & &\ddots\\
&            &       &\aa_0       &\dots  &\dots &\dots&\dots  &\aa_k\\
\bf0&&&\aa_0^{(\nu+1)}     &\dots  &\dots &\dots&\dots  &\aa_k^{(\nu+1)}\\
&&&\vdots      &       &      &       &&\vdots\\
&&&\aa_0^{(k)}       &\dots  &\dots &\dots&\dots  &\aa_k^{(k)}
\endpmatrix,\end{equation}

\no it turns out that the discrete problem, at a given point $t_n=nh$, can
be cast in matrix form as
\begin{equation}\label{discr}(\calA_2 \otimes I_m) \bfy_{new}  = h\bff_{new}
-(\calA_1\otimes I_m)\bfy_{old}.\end{equation}

\no This corresponds to a GLM of the form (\ref{glm}) with
coinciding internal and external stages, and $A = \calA_2^{-1}$,
$U = -\calA_2^{-1}\calA_1$ (indeed, matrix $\calA_2$ turns out to
be nonsingular). Such methods were called {\em Block BVMs with
Memory (B$_2$VM$_2$s)} in \cite{BrTr98}. In the GLM notation, the
abscissae defining this GLM are $c_i = i$, $i=1,2,\dots,r$.
However, in order to guarantee the $L$-stability of the
corresponding method, it is sometimes required to change them into
the following set of abscissae:
\begin{eqnarray}\nonumber
c_i &=& i, \qquad i=1,\dots,\ell-1, \\  \label{intsteps}
c_{\ell+j} &=&
\ell-1+\sum_{m=0}^j\xi_m, \qquad j=0,\dots,r-\ell,\\
\nonumber
&&\mbox{with positive}\quad \{\xi_m\}\quad\mbox{s.t.}\quad c_r
=\ell.\end{eqnarray}

A first possible choice, as also suggested in \cite{BrTr98}, is
that of setting \begin{equation}\label{choice1} \xi_m =
\zeta^{m+1}, \qquad m = 0,\dots,r-\ell,\end{equation}

\no where $\zeta$ is the positive root of the polynomial
$$
p(z) = \sum_{m=0}^{r-\ell}z^{m+1} -1.
$$

\no However, in general such a value of $\zeta$ is an irrational
number, which implies that the coefficients of the corresponding
GLM are irrational numbers as well. In order to have methods whose
coefficients are always rational numbers, a slightly different
choice for the abscissae can be made, i.e.:
\begin{equation}\label{choice2}
\xi_m = \frac{ 2^{r-\ell-m}}{2^{r-\ell+1}-1}, \qquad m =
0,\dots,r-\ell.
\end{equation}

Clearly, when ~$\ell=r$~ we have a uniform mesh, whereas for
~$\ell<r$~ the last ~$r-\ell+1$~ stepsizes are geometrically
decreasing, in the case of choice (\ref{choice1}), or
approximately halved, in the case of choice (\ref{choice2}). In
such a case, the points corresponding to the abscissae
~$\{c_0,\dots,c_{\ell-1},c_r\}$~ are equally spaced and will be
the ones needed for the subsequent integration step. This implies
that $r-\ell$ zero columns must be appropriately inserted in
matrix $\calA_1$ in (\ref{A1A2}), i.e., those referring to the
approximations at the abscissae $\{c_\ell,\dots,c_{r-1}\}$
(according to \cite{BrTr98}, such points are called {\em auxiliary
points}). Consequently, such a GLM based on GBDF is uniquely
determined by the triple of integers $(k,r,\ell)$. Clearly, the
order of the corresponding method is $p=k$. In addition, for all
practical values of $k$, by appropriately choosing the values of
$(r,\ell)$, it is possible to guarantee the $L$-stability of such
methods. However, it turns out that, for each value of $k$, the
couples $(r,\ell)$ are not unique. Consequently, an additional
criterion will be considered, in the next section, which is aimed
to speed-up the iterative solution of the corresponding discrete
problems, via the {\em blended implementation} of the methods.

\section{Blended General Linear Methods}\label{blend}

In the previous section, we devised a procedure for obtaining a
whole class of implicit GLMs based on GBDF, containing $L$-stable
methods of arbitrarily high order. We now consider the problem of
efficiently solving the generated discrete problems which, at a
given time step, assume the form (see (\ref{glm})),
\begin{equation}\label{simdis}
\bfy -h(A\otimes I_m)\bff = \bfeta,
\end{equation}

\no where, according to (\ref{A1A2}) and (\ref{discr}),
$$\bfy = \bfy_{new}, \qquad \bff = \bff_{new}, \qquad
A = \calA_2^{-1}, \qquad \bfeta=-(\calA_2^{-1}\calA_1\otimes I_m) \bfy_{old}.$$

\no A discrete problem in the form (\ref{simdis}) can be
efficiently solved, under suitable hypotheses, via the {\em
blended implementation} of the method. We recall that the blended
implementation of block implicit methods has been previously
considered in the framework of one-step methods
\cite{Br00,BrMa02,BrMa04,BrMa07,BrMa08,BrMaMu06,BrTr01} and has
been implemented in the computational codes {\tt BiM}
\cite{BrMa04,bim} and {\tt BiMD} \cite{BrMaMu06,bim}. In order to
describe the {\em blended implementation} of the GLM
(\ref{simdis}), it is convenient to consider its application for
solving the usual test equation (\ref{test}). In such a case, in
fact, the discrete problem reduces to a linear system of dimension
$r$:
\begin{equation}\label{prob} (I-qA)\bfy = \bfeta, \qquad
q=h\lam,\end{equation}

\no where, hereafter, $I$ denotes the identity matrix of dimension $r$. Such
a linear system is clearly equivalent
to\begin{equation}\label{prob1} \cc(A^{-1} -qI)\bfy = \cc
A^{-1}\bfeta \equiv \bfeta_1, \end{equation}

\no with $\cc>0$ a free parameter. By introducing the {\em weight
function} \begin{equation}\label{teta} \ttt(q) = (1-\cc q)^{-1}
I,\end{equation} we can then obtain an equivalent discrete problem by
combining, with weights $\ttt(q)$ and $I-\ttt(q)$, respectively, the linear
systems
(\ref{prob})-(\ref{prob1}):
\begin{eqnarray}\nonumber
M(q)\bfy &\equiv& \left(\ttt(q)(I-qA) +\ \cc
(I-\ttt(q))(A^{-1}-qI) \right)\bfy \\ &=& \ttt(q)\bfeta +
(I-\ttt(q))\bfeta_1 ~\equiv~ \bfeta(q).\label{Mq}\end{eqnarray}

\no Equation (\ref{Mq}) defines the {\em Blended General Linear
Method (Blended GLM)} corresponding to the original GLM (\ref{prob}). Then,
in view of the fact that $$M(q)\approx \left\{ \begin{array}{cl}I,
&\mbox{for}\quad q\approx0,\\[2mm]
-\cc q I, &\mbox{for}\quad |q|\gg1,\end{array}\right.$$

\no the following {\em blended iteration} for solving the discrete problem
(\ref{Mq}) is naturally induced:
\begin{equation}\label{Nq}
N(q)\bfy^{(i+1)} \equiv (I-\cc q I)\bfy^{(i+1)} =
(N(q)-M(q))\bfy^{(i)} +\bfeta(q), \qquad
i=0,1,\dots.\end{equation}

\begin{rem} By considering that in (\ref{Nq}) $N(q)^{-1}=\ttt(q)$ (see
(\ref{teta})), it is quite straightforward to realize that, in the case of
problem (\ref{ivp}), the corresponding blended iteration formally becomes
\begin{eqnarray}\nonumber
\bfdel^{(i)} &=& N^{-1}\left( \ttt\left( (I-\cc A^{-1})\otimes I_m\,
(\bfy^{(i)}-\bfeta) -h(A-\cc I)\otimes I_m\bff^{(i)}\right)\right.\\
\label{blendgen}&&\left.+\cc\left(A^{-1}\otimes I_m\,
(\bfy^{(i)}-\bfeta)-hI\otimes I_m\,\bff^{(i)}\right)\right),\\
\nonumber \bfy^{(i+1)} &=& \bfy^{(i)}-\bfdel^{(i)}, \qquad
i=0,1,\dots,
\end{eqnarray}

\no where \begin{equation}\label{N}
N \equiv \ttt^{-1} = I\otimes (I_m-h\cc J),\end{equation}

\no with $J$ the Jacobian of $f$ evaluated at the last known point.
Consequently, only the factorization of {\em one matrix of the size of the
continuous problem} is required.
\end{rem}

Coming back to the iteration (\ref{Nq}), we observe that (see
\cite{BrMa02,BrMa08}) the iteration matrix corresponding to the
blended iteration (\ref{Nq}) turns out to be given by
\begin{equation}\label{Zq} Z(q) \equiv I-N(q)^{-1}M(q) =
\frac{q}{(1-\cc q)^2} A^{-1}(A-\cc I)^2.\end{equation}

\no According to the linear analysis described in \cite{BrMa08},
we shall now study the spectral properties of $Z(q)$, in order to
obtain a linear analysis of convergence for the iteration
(\ref{Nq}). For this purpose, hereafter let $\rho(q)$ denote the
spectral radius of $Z(q)$. Clearly, the iteration will be
convergent if and only if $\rho(q)<1$; moreover, the set $$\Gamma
= \left\{ q\in\CC \,:\, \rho(q)<1\right\}$$

\no is the {\em region of convergence} of the iteration. For the
sake of completeness (see  \cite{BrMa08} for details), we recall
that the iteration is:
\begin{itemize}

\item {\em $0$-convergent}, if $\rho(0)=0$;

\item {\em $A$-convergent}, if $\CC^-\subseteq\Gamma$;

\item {\em $L$-convergent} with {\em index} $\nu_\infty$,
if it is $A$-convergent,
$$Z_\infty\equiv \lim_{q\rightarrow\infty}Z(q) = O,$$ and
$\nu_\infty$ is the index of nilpotency of $Z_\infty$.

\end{itemize}
Clearly, the iteration (\ref{Nq}) is {\em $0$-convergent}, since
$Z(0)=O$. Moreover, one easily verifies (from (\ref{Zq})) that
$$Z(q)\rightarrow O, \quad\mbox{as}\quad q\rightarrow\infty,$$

\no so that $A$-convergence and $L$-convergence (with index 1) are
equivalent. In addition to this,
$$\rho(q) \approx\left\{
\begin{array}{ll}
\approx \tro q, &\mbox{for}\quad  q\approx0,\\[2mm]
\approx \troinf q^{-1}, &\mbox{for}\quad  |q|\gg1,
\end{array}\right.$$

\no where $\tro$ and $\troinf$ are the {\em nonstiff amplification
factor} and the {\em stiff convergence factor} of the iteration,
respectively. Finally, for a $0$-convergent iteration,
$A$-convergence is equivalent to requiring that the {\em maximum
amplification factor}, $$\rostar = \max_{x>0}\rho(ix),$$

\no (with $i$ denoting, as usual, the imaginary unit) is not
greater than 1. According to \cite{BrMa08}, $\tro$, $\troinf$, and
$\rostar$ are evaluation parameters for the blended iteration
(\ref{Nq}) (the smaller they are, the better the properties of the
iteration). Consequently, if $\rostar\le1$, the iteration is
$L$-convergent and, therefore, appropriate for $L$-stable methods
\cite{BrMa08}, like the GLMs defined in Section~\ref{gbdf}.

Because of the form (\ref{Zq}) of the iteration matrix, it
turns out that (see also \cite{BrMa08}):
$$\tro = \rho\left( A^{-1}(A-\cc I)^2 \right), \qquad \troinf =
\cc^{-2}\tro, \qquad \rostar = (2\cc)^{-1} \tro,$$

\no with $\rho\left( A^{-1}(A-\cc I)^2 \right)$ denoting the
spectral radius of that matrix. Consequently, the value of the free positive
parameter $\cc$ is chosen in order to minimize $\rostar$. In
Table~\ref{tab1}, the relevant figures are reported for selected GLMs based
on GBDF, each characterized by the corresponding triple $(k,r,\ell)$, when
the choice (\ref{choice1}) for the auxiliary points is considered.
Similarly, in Table~\ref{tab2}, the same results, obtained with
the choice (\ref{choice2}) of the auxiliary points, are listed. In both
tables, we list $r-\ell$ (i.e., the number of the auxiliary points)
in place of $\ell$. The specified values of the blocksize $r$, have been
here chosen as small as possible, in order to minimize the computational
cost per step.

As it can be seen, the blended iteration corresponding
to each method turns out to be $L$-convergent. Moreover, the value of the
parameter $\cc$ (with the only exception of that corresponding to
$k=4$, for both the choices of the auxiliary steps) turns out to
coincide with the value $$\cc^* = \min_{\mu\in\sigma(A)} |\mu|,$$

\no suggested in \cite{BrMa02} for a class of Blended Implicit
Methods (see also \cite{BrMa08}).

As an example, we list the matrices $A$ and $U$ which define the
third order GLM corresponding to the triple $(3,2,2)$, requiring no
auxiliary steps (hereafter, let $\bfc\in\RR^r$ denote the vector of the
abscissae),
\begin{equation}\label{322}
\bfc = (\,1,\,2\,)^T,\qquad A = \frac{1}{23}\pmatrix{rr} 22 &-4\\
36 &6\endpmatrix, \qquad U = \frac{1}{23}\pmatrix{rr} -5 &28\\ -4
&27\endpmatrix,
\end{equation}

\no and those defining the fourth order GLM corresponding to the
triple $(4,4,3)$, with the choice (\ref{choice1}),
\begin{eqnarray*}
&&\bfc = \left(\, 1,\, 2,\, \frac{3+\sqrt{5}}{2},\,
3\,\right)^T,\\~\\ &&A =\nonumber\\[2mm] &&\hspace{-.25cm}{\small
\pmatrix{rrrr}
0.85795933248329& -0.22588594615620&  0.09292560797716& -0.02384741384935\\
1.14013555838905&  0.75295315385401& -0.30975202659054&  0.07949137949784\\
1.16454145194449&  0.91413597782316&  0.27338586108581& -0.07015876367984\\
1.16304178987375&  0.90057211551936&  0.50490826434742&  0.09507592794253
\endpmatrix},\nonumber\\~\\
&&U = {\small\pmatrix{rrrr}
0.07149661104027& -0.44184164162566&  0&  1.37034503058538\\
0.09501129653242& -0.52719452791448&  0&  1.43218323138206\\
0.09704512099537& -0.53021970356702&  0&  1.43317458257165\\
0.09692014915615& -0.53024220062923&  0&  1.43332205147308
\endpmatrix}.\nonumber\end{eqnarray*}
\no and with the choice (\ref{choice2}),

\begin{eqnarray*}
\bfc &=& \left(\, 1,\, 2,\, \frac{8}{3},\,3\,\right)^T,\\~\\ A &=&
\frac{1}{6336684}\, \pmatrix{rrrr}
5429268 & -1381941 &   570807 &  -174960\\
7249176 &  4606470 & -1902690 &   583200\\
7421568 &  5690784 &  2388672 &  -732160\\
7415388 &  5637357 &  3628233 &   198936
\endpmatrix,\nonumber\\~\\
U &=& \frac{1}{6336684}\, \pmatrix{rrrr}
452439   &  -2798388  &  0  &  8682633 \\
604098   &  -3345408  &  0  &  9077994 \\
618464   &  -3365888  &  0  &  9084108 \\
617949   &  -3366036  &  0  &  9084771
\endpmatrix,\nonumber\end{eqnarray*}

\no which slightly differ from each other.

\begin{table}[t]\centerline{
\begin{tabular}{|r|r|c|c|c|c|c|}
\hline
$k$& $r$&  $r-\ell$& $\cc$& $\tro$& $\troinf$& $\rostar$\\
\hline
  3&  2&  0&   0.7223  &   0.2272  &   0.4355  &   0.1573 \\
  4&  4&  1&   0.6195  &   0.3802  &   0.9908  &   0.3069 \\
  6&  5&  1&   0.6063  &   0.5734  &   1.5600  &   0.4729 \\
  8&  6&  1&   0.5769  &   0.6380  &   1.9170  &   0.5530 \\
 10&  7&  1&   0.5502  &   0.6626  &   2.1887  &   0.6021 \\
 12&  9&  2&   0.5271  &   0.7345  &   2.6438  &   0.6968 \\
 14& 10&  2&   0.5127  &   0.7366  &   2.8022  &   0.7183 \\
 16& 11&  2&   0.4999  &   0.7345  &   2.9393  &   0.7347 \\
\hline
\end{tabular}}
\caption{Parameters of the blended iteration associated with Blended GLMs,
based on GBDF, characterized by the triple $(k,r,\ell)$ and the choice
(\ref{choice1}).}
\label{tab1}
\bigskip
\centerline{
\begin{tabular}{|r|r|c|c|c|c|c|}
\hline
$k$& $r$&  $r-\ell$& $\cc$& $\tro$& $\troinf$& $\rostar$\\
\hline
  3&  2&  0&  0.7223  &  0.2272  &  0.4355  &  0.1573 \\
  4&  4&  1&  0.6249  &  0.3827  &  0.9801  &  0.3062 \\
  6&  5&  1&  0.6082  &  0.5740  &  1.5520  &  0.4719 \\
  8&  6&  1&  0.5778  &  0.6381  &  1.9113  &  0.5522 \\
 10&  7&  1&  0.5507  &  0.6625  &  2.1845  &  0.6015 \\
 12&  9&  2&  0.5274  &  0.7345  &  2.6407  &  0.6964 \\
 14& 10&  2&  0.5130  &  0.7366  &  2.7998  &  0.7180 \\
 16& 11&  2&  0.5000  &  0.7345  &  2.9374  &  0.7344 \\
\hline
\end{tabular}}
\caption{Parameters of the blended iteration associated with Blended GLMs,
based on GBDF, characterized by the triple $(k,r,\ell)$ and the choice
(\ref{choice2}).}
\label{tab2}
\end{table}

\begin{figure}[ht]
\begin{center}
\includegraphics*[width=10cm]{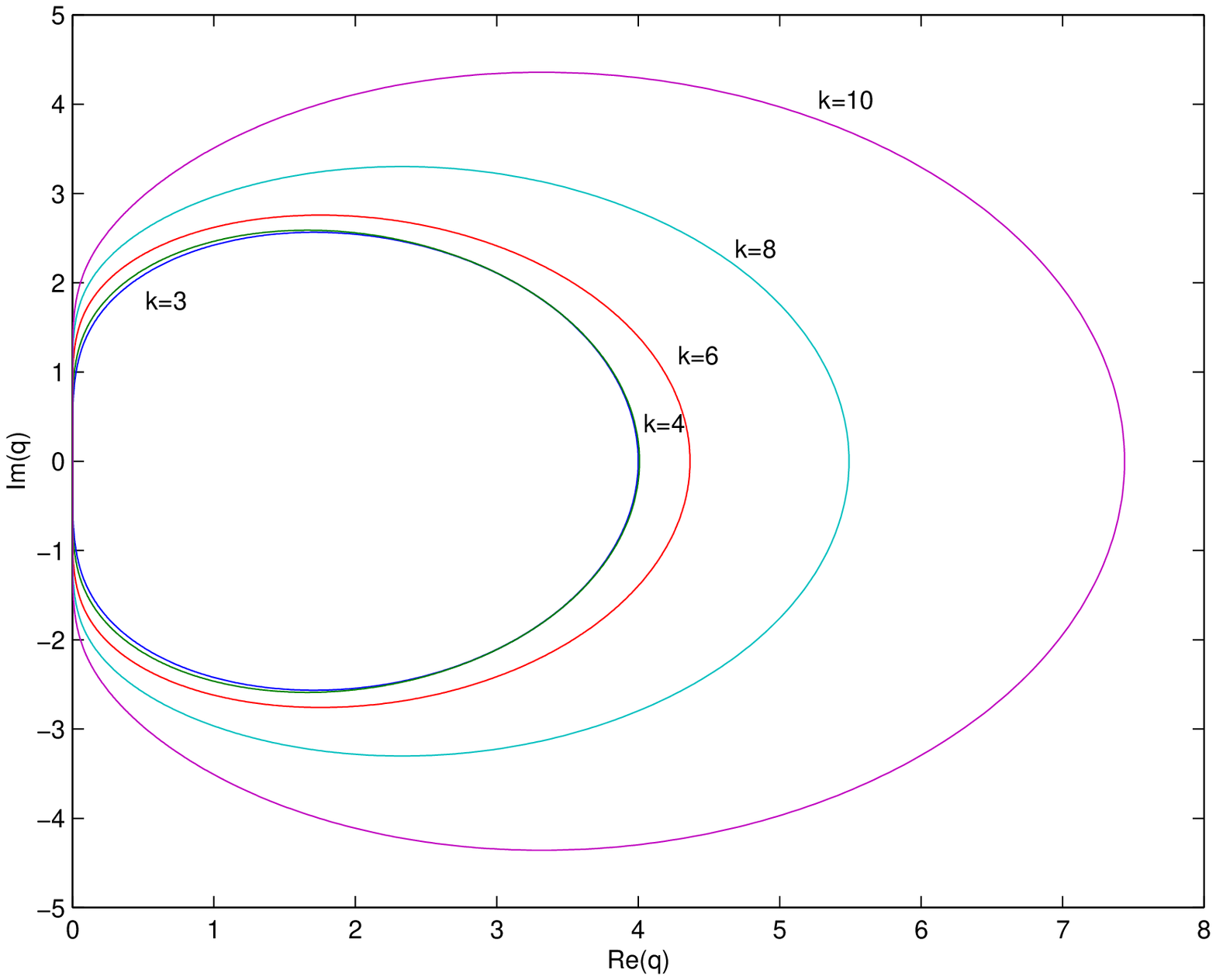}
\caption{Boundary loci of the GLMs obtained with the choice
(\ref{choice1}), $k=3,4,6,8,10$.} \label{fig1}
\bigskip
\includegraphics*[width=10cm]{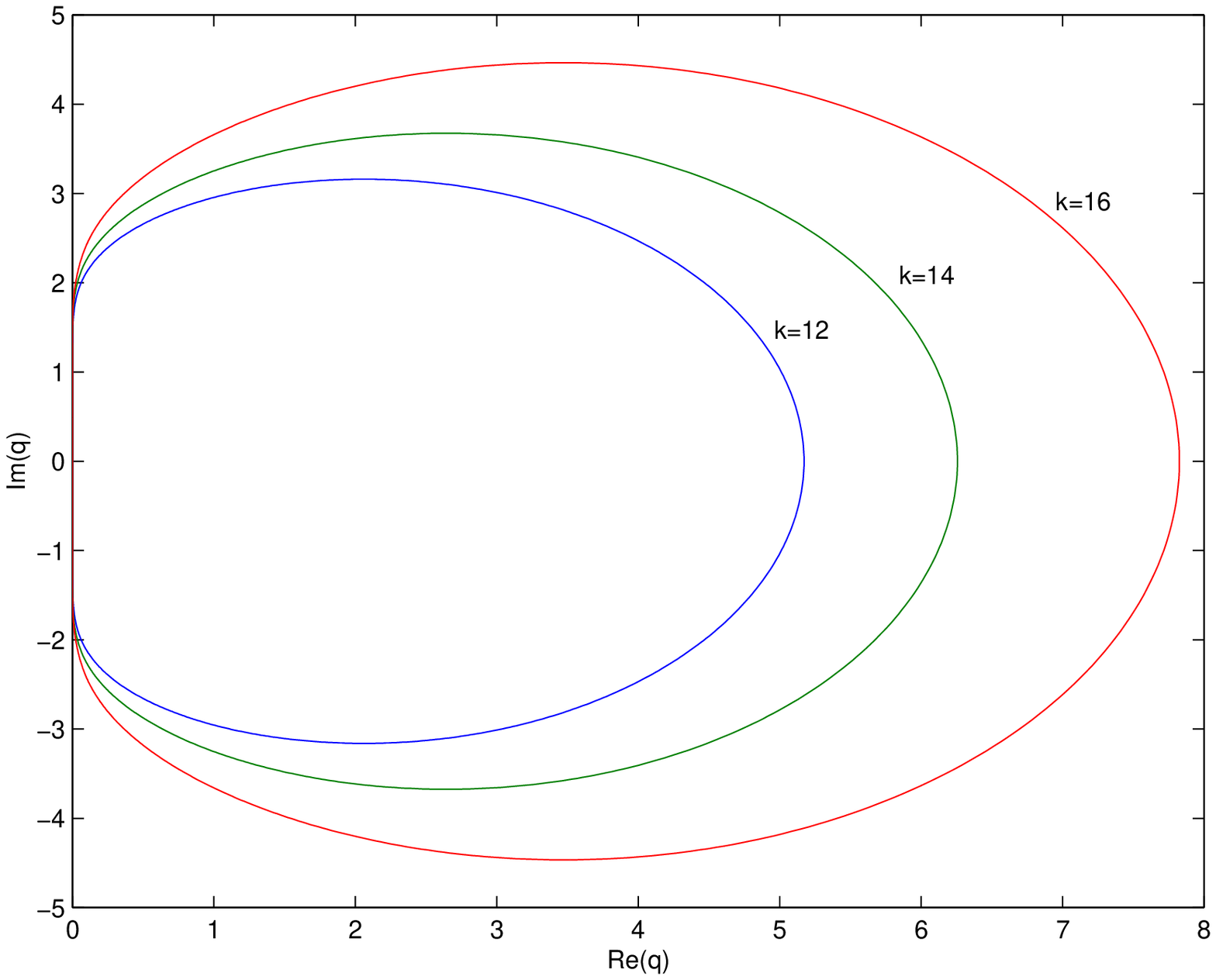}
\caption{Boundary loci of the GLMs obtained with the choice
(\ref{choice1}), $k=12,14,16$.}
 \label{fig2}
\end{center}
\end{figure}

\begin{figure}[ht]
\begin{center}
\includegraphics*[width=10cm]{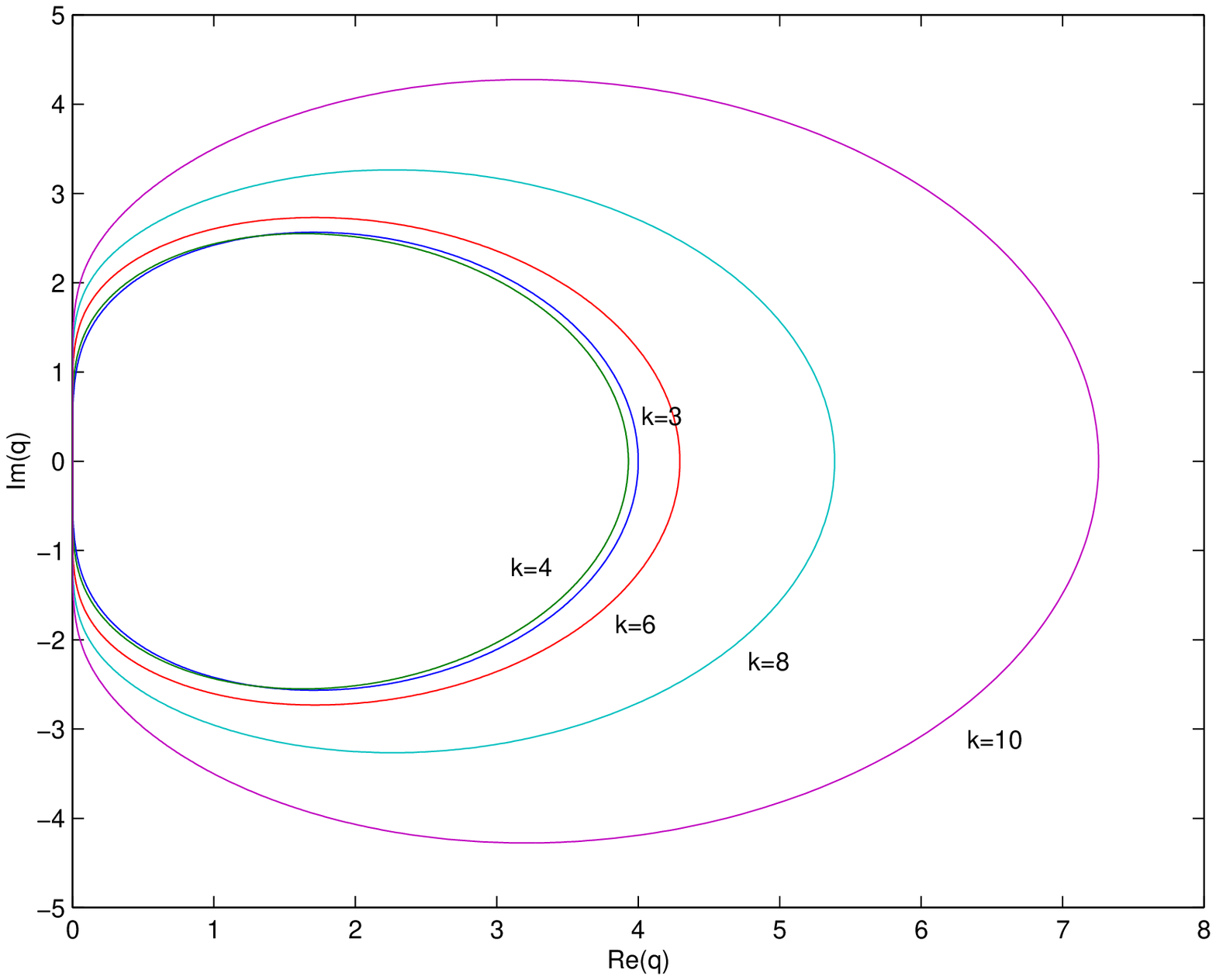}
\caption{Boundary loci of the GLMs obtained with the choice
(\ref{choice2}), $k=3,4,6,8,10$.} \label{fig3}
\bigskip
\includegraphics*[width=10cm]{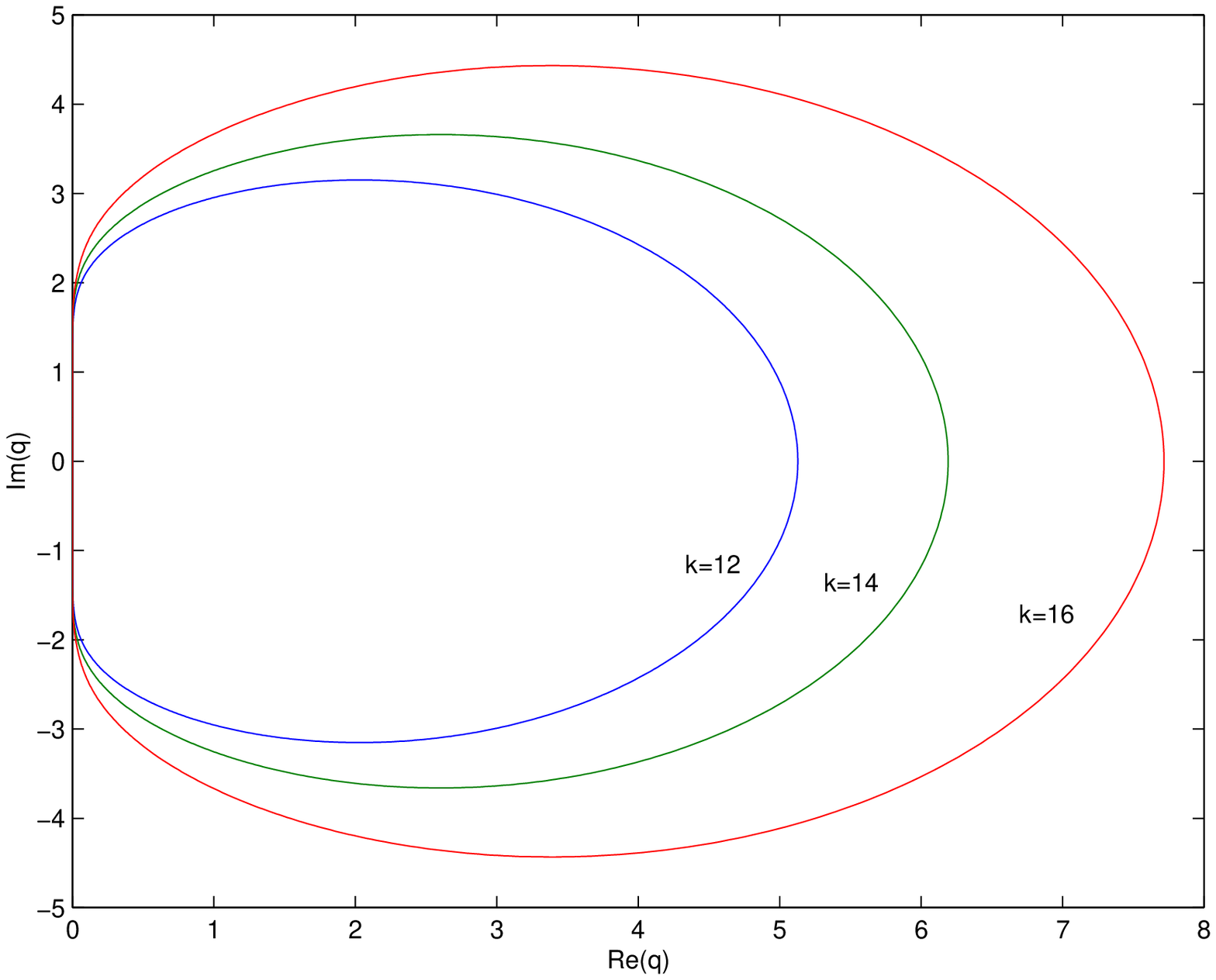}
\caption{Boundary loci of the GLMs obtained with the choice
(\ref{choice2}), $k=12,14,16$.}
 \label{fig4}
 \end{center}
\end{figure}

Finally, in Figures~\ref{fig1} and \ref{fig2}, we plot the
boundary of the stability regions of the methods listed in
Table~\ref{tab1}, whereas in Figures~\ref{fig3} and \ref{fig4}, we
plot the boundary of the stability regions of the methods listed
in Table~\ref{tab2}.

\begin{figure}[t]
\begin{center}
\includegraphics*[width=11cm]{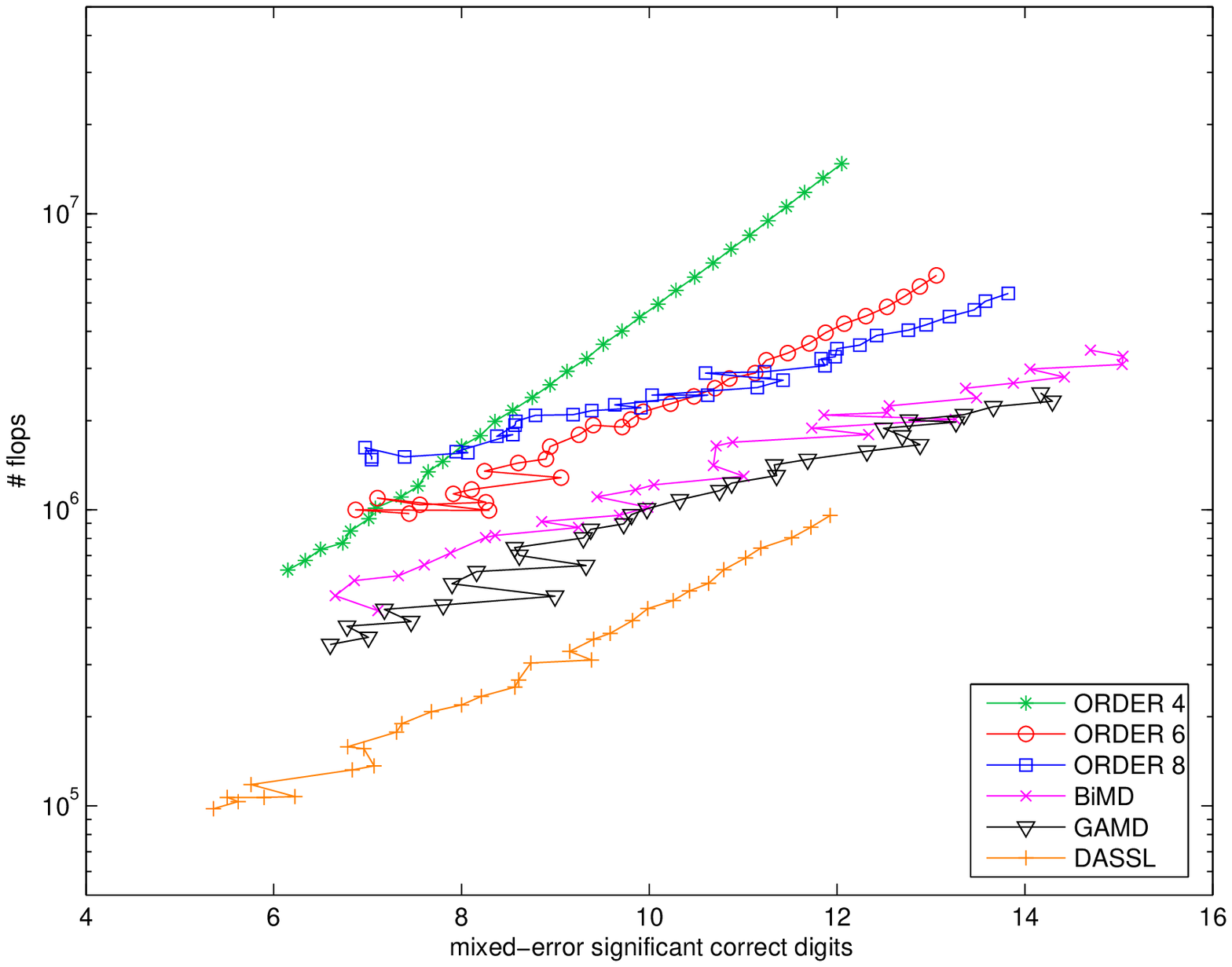}
\caption{Numerical results for the Pollution Problem.}
\label{fig5}
\bigskip
\includegraphics*[width=11cm]{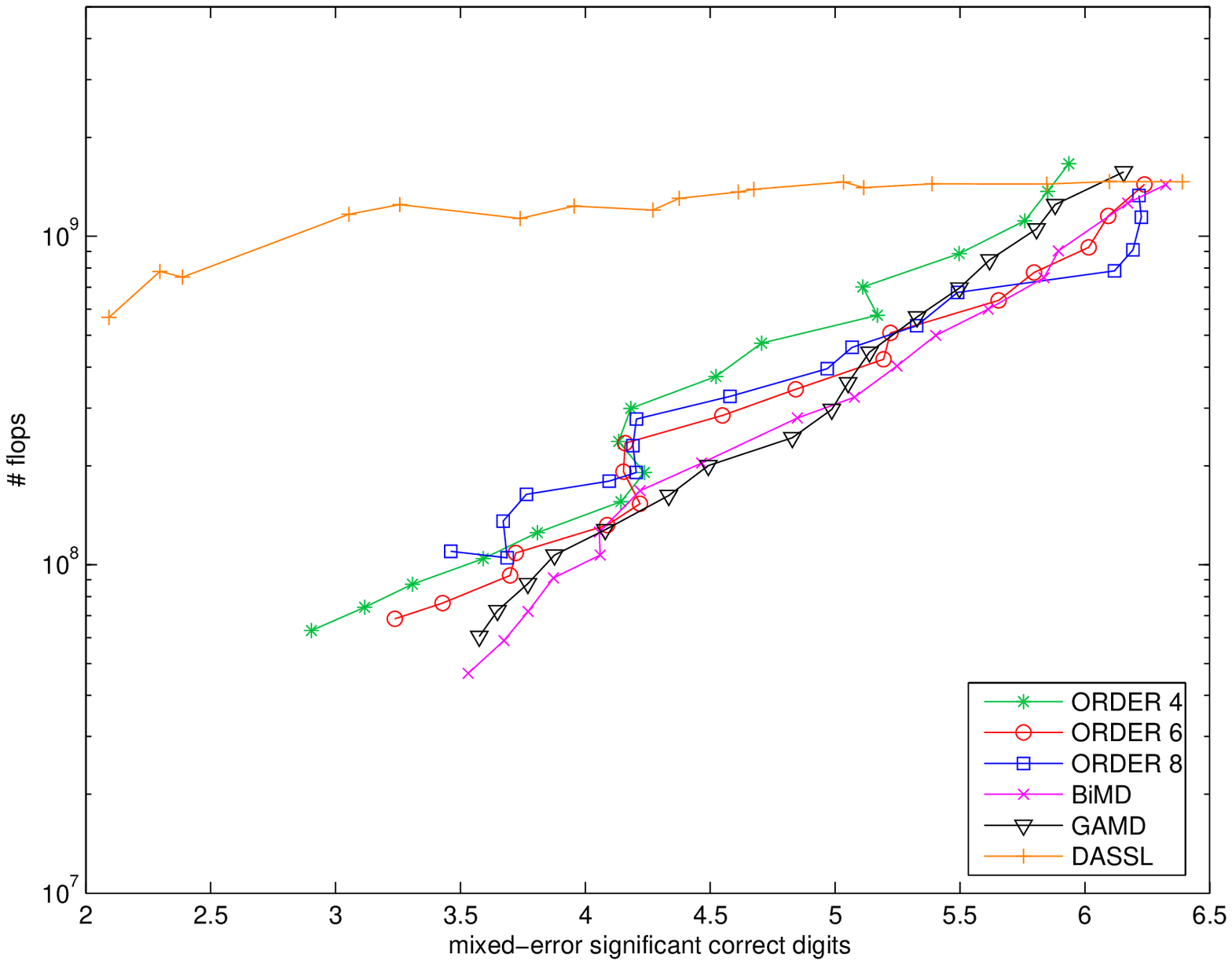}
\caption{Numerical results for the Elastic Beam Problem.}
\label{fig6}
\end{center}
\end{figure}

\begin{figure}[t]
\begin{center}
\includegraphics*[width=11cm]{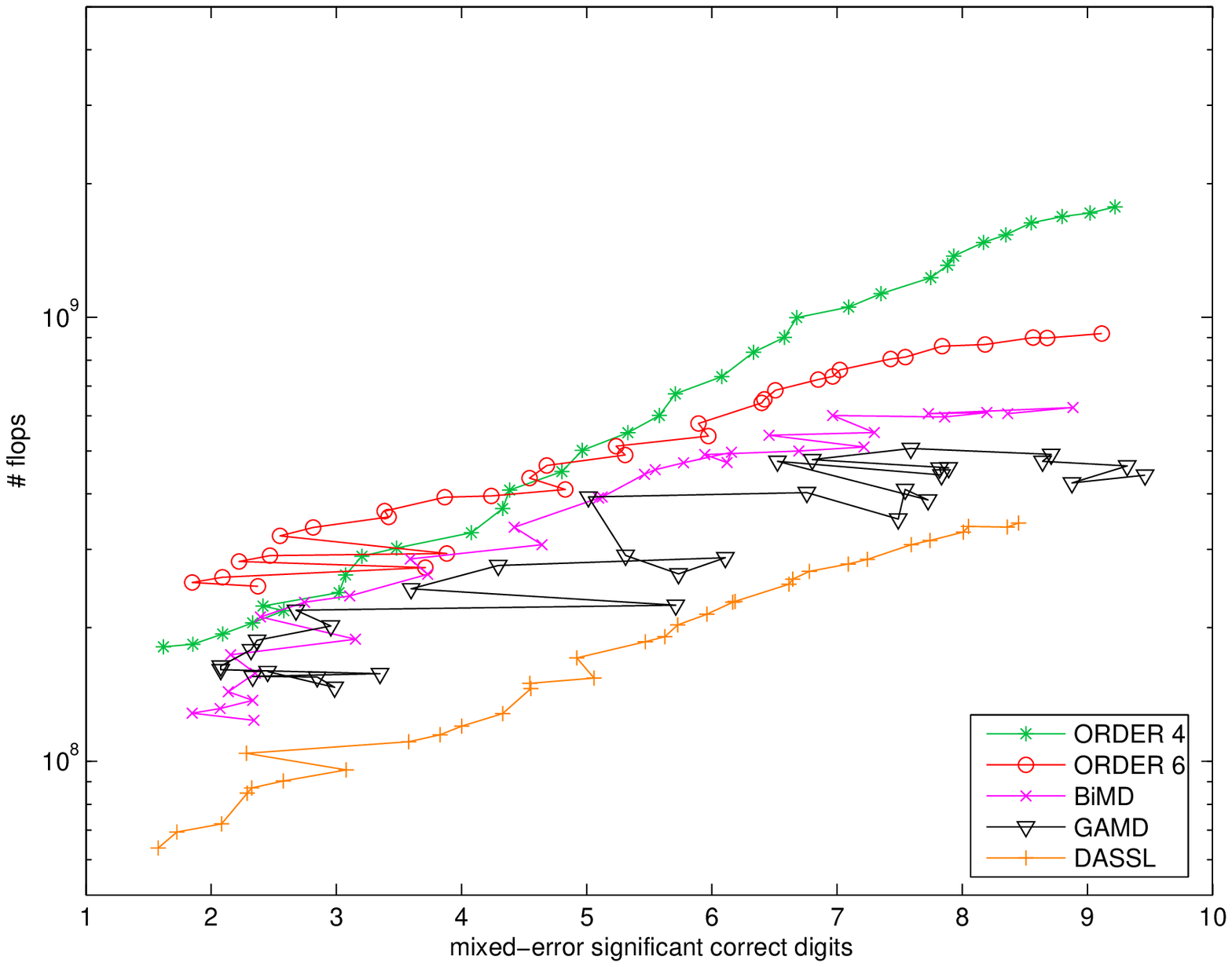}
\caption{Numerical results for the Emep Problem.}\label{fig7}
\bigskip
\includegraphics*[width=11cm]{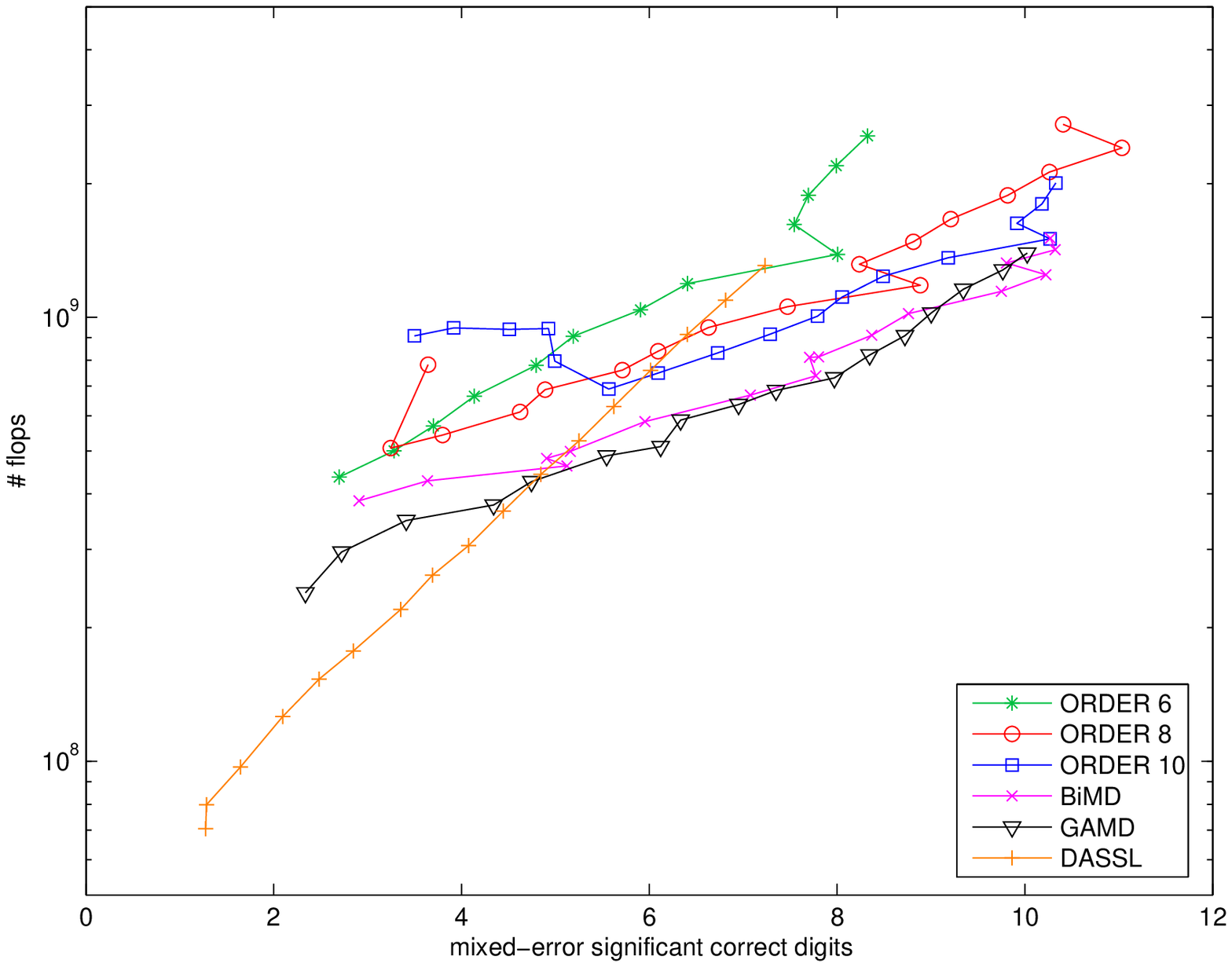}
\caption{Numerical results for the Ring Modulator
Problem.}\label{fig8}
\end{center}
\end{figure}

\section{Implementation details}\label{impl}

In the actual implementation of the above Blended GLMs, three main points
need to be clarified:
\begin{itemize}
 \item the choice of a suitable starting procedure;
 \item an efficient local error estimate;
 \item the variable stepsize implementation of the methods.
\end{itemize}
All of them are briefly sketched here.

Concerning the first issue, namely the definition of a starting
procedure to obtain the first vector of approximations, from the
initial condition $y_0$ of problem (\ref{ivp}), a natural
candidate is the block GBDF of order $k$ and blocksize $r=k$. In
more detail, at the very beginning, one solves the discrete
problem
\begin{equation}\label{sp}\calA \otimes I_m \bfy -h{\cal B}\otimes I_m
\bff = \bf0,\end{equation}

\no where$$\bfy = (~y_0,~y_1,~\dots,~y_k~)^T, \qquad
\bff =(~f_0,~f_1,~\dots,~f_k~)^T,$$

\no and, by denoting with $\{\aa_i^{(j)}\}$ the $\aa$-s
coefficients of the additional initial and final methods
(\ref{init}) and (\ref{fine}), the matrices $\calA$ and ${\cal B}$
are defined as (see (\ref{AA})):

$$\calA = \pmatrix{l|lcl}
\aa_0^{(1)} & \aa_1^{(1)} &\dots     &\aa_k^{(1)}\\
\vdots      & \vdots        & &\vdots\\
\aa_0^{(\nu-1)} &\aa_1^{(\nu-1)} &\dots     &\aa_k^{(\nu-1)}\\
\aa_0  &\aa_1     &\dots     &\aa_k\\
\aa_0^{(\nu+1)} &\aa_1^{(\nu+1)} &\dots &\aa_k^{(\nu+1)}\\
\vdots      & \vdots &      &\vdots\\
\aa_0^{(k)} & \aa_1^{(k)} &\dots &\aa_k^{(k)}
\endpmatrix, ~ {\cal B} = (\,{\bf0} \,|\,
I_k\,) ~ \in~\RR^{k\times k+1}.$$

\no For the efficient solution of equation (\ref{sp}), a corresponding
blended methods, with the associated blended iteration, can be conveniently
and easily defined.

The estimate of the local error is obtained by {\em deferred
correction}, as described in \cite[Chapter~10]{BrTr98} (see also
\cite{BrMa04,BrMa06}), by first ``plugging'' the local discrete
solution (see (\ref{discr})) $(\,\bfy_{old},\,\bfy_{new}\,)^T$ in
the discrete problem defined by a higher order method. In such a
way, one obtains an estimate $\bftau$ of the local truncation
error. The same is done by considering an equivalent formulation
of the same method (compare with (\ref{simdis}) and
(\ref{prob})-(\ref{prob1})), thus obtaining a new approximation
$\bftau_1=\cc A^{-1}\otimes I_m \bftau$. After that, the estimate $\bfe$ of
the local error is obtained by performing one blended iteration, namely by
formally solving the block diagonal linear system (see (\ref{N}))
$$N\bfe = \theta \bftau +(I-\theta)\bftau_1.$$

\no In more detail, when using the GLM defined by the triple
$(k,r,\ell)$, the corresponding method used for estimating
$\bftau$ is the GBDF defined by the triple $(k+1,r,\ell)$.
Consequently, the cost for estimating $\bfe$ is the same as that
for carrying out one blended iteration. This can be done for all
methods listed in Tables~\ref{tab1} and \ref{tab2}, with the only
exception being the triple (3,2,2), for which the value of $r$
should be increased. For this reason, such a method has not been
considered for carrying out the numerical tests in
Section~\ref{numer}.

Finally, for the efficient variation of the stepsize, we have
considered the Nordsiek implementation of the methods, firstly
introduced in \cite{Nord} and later modified according to, e.g.,
\cite{La93,RaHi93}. For the sake of completeness, we also mention
that the initial guess for the blended iteration
(\ref{blendgen})-(\ref{N}) is obtained by extrapolation, through
the interpolating polynomial on (see (\ref{discr})) $\bfy_{old}$.

\section{Numerical Tests}\label{numer} In this section, we report a few
numerical tests comparing a Matlab fixed-order implementation of
the Blended GLMs presented here, with some of the most reliable
codes currently available, on selected stiff test problems. Both
the problems and the codes have been taken from the current
release (release~2.4) of Test Set for IVP Solvers \cite{TestSet}.
In particular, we have considered the following solvers:
\begin{itemize}
 \item BiMD, which uses a blended iteration for solving the generated
discrete problems;
\item GAMD, which uses a nonlinear splitting for solving the
discrete problems, generated by block BVMs in the Generalized
Adams Methods (GAMs) family;
\item DASSL, which implements standard BDF methods.
\end{itemize}

\no The problems considered are:
\begin{itemize}
 \item Pollution (of dimension 20);
 \item Elastic Beam (of dimension 80);
 \item Emep (of dimension 66);
 \item Ring Modulator (of dimension 15).
\end{itemize}

\no The comparisons have been summarized in corresponding {\em
work-precision diagrams} \cite{TestSet}, where the computational
cost is plotted versus accuracy. A standardized cost has been
computed as the number of floating-point operations for the
factorizations and the system solvings required by each code,
while accuracy has been measured in terms of {\em mixed-error
significant correct digits} \cite{TestSet}. For all codes, the
used tolerances are essentially those specified in the Test Set.
Figures~\ref{fig5}--\ref{fig8} summarize the results obtained,
where the label ``ORDER $k$'' is used for the fixed-order
implementation of the Blended GBDF of order $k$ (see
Table~\ref{tab2}): as one can see, for each problem, the selected
fixed-order implementations of the Blended GBDF presented here,
appear to be competitive with the above mentioned solvers.

\section{Conclusions}\label{concl} In this paper, a straightforward
approach for deriving $L$-stable General Linear Methods (GLMs) of
arbitrarily high order has been presented. Such an approach relies on the
class of Boundary Value Methods (BVMs) for ODEs. In the same framework
corresponding starting procedures are easily derived, as well as appropriate
error estimates.

The generated discrete problems can be efficiently solved by means
of the blended implementation of the methods, thus defining
corresponding Blended GLMs, equivalent to the original methods,
from the point of view of the stability and accuracy
properties. The corresponding blended iterations are all
$L$-convergent, thus appropriate for the underlying $L$-stable
methods.

A number of numerical tests, on problems taken from the Test Set for IVP
Solvers, prove that the obtained methods are competitive with some of the
best codes currently available.

The availability of methods having arbitrarily high-order makes
them good candidates for an efficient variable-order
implementation.


\begin{thebibliography}{99}

\bibitem{AcTr02} L.\,Aceto, D.\,Trigiante. On the $A$-stable methods in the
GBDF class. {\em Nonlinear Anal. Real World Appl.} {\bf 3} (2002) 9--23.

\bibitem{Br00}
L.\,Brugnano. Blended Block BVMs (B$_3$VMs): a family of
economical implicit methods for ODEs. {\em Jour. Comput. Appl.
Mathematics} {\bf 116} (2000) 41--62.

\bibitem{BrMa02}
L.\,Brugnano, C.\,Magherini. Blended Implementation of Block Implicit
Methods for ODEs. {\em Appl. Numer. Math.}  {\bf 42} (2002) 29--45.

\bibitem{BrMa04}
L.\,Brugnano, C.\,Magherini. The BiM code for the numerical
solution of ODEs. {\em Jour. Comput. Appl. Mathematics}  {\bf
164--165} (2004) 145--158.

\bibitem{BrMa06}
L.\,Brugnano, C.\,Magherini. Economical Error Estimates for Block Implicit
Methods for ODEs via Deferred Correction. {\em Appl. Numer. Math.}  {\bf 56}
(2006) 608--617.

\bibitem{BrMa07}
L.\,Brugnano, C.\,Magherini. Blended Implicit Methods for solving
ODE and DAE problems, and their extension for second order
problems. {\em Jour. Comput. Appl. Mathematics} {\bf 205} (2007)
777--790

\bibitem{BrMa08}
L.\,Brugnano, C.\,Magherini. Recent advances in linear analysis of
convergence for splittings for solving ODE problems. {\em Appl.
Numer. Math.} {\bf 59} (2009) 542--557.

\bibitem{BrMa08-1}
L.\,Brugnano, C.\,Magherini.  Blended General Linear Methods based on
Generalized BDF. {\em AIP Conference Proceedings} {\bf 1048} (2008)
871--874.

\bibitem{BrMaMu06}
L.\,Brugnano, C.\,Magherini, F.\,Mugnai. Blended Implicit Methods
for the numerical solution of DAE problems. {\em Jour. Comput.
Appl. Mathematics}  {\bf 189}  (2006) 34--50.

\bibitem{BrTr96}
L.\,Brugnano, D.\,Trigiante. Convergence and stability of Boundary
Value Methods for ordinary differential equations. {\em Jour.
Comput. Appl. Mathematics} {\bf 66} (1996) 97--109.

\bibitem{BrTr98}
L.Brugnano, D.Trigiante. {\em Solving Differential Problems by
Multistep Initial and Boundary Value Methods}, Gordon and Breach
Science Publ., 1998.

\bibitem{BrTr01}
L.\,Brugnano, D.\,Trigiante. Block Implicit Methods for ODEs, in
{\em Recent Trends in Numerical Analysis}, D.Trigiante ed., Nova
Science Publ. Inc., New York, 2001, pp. 81-105.

\bibitem{BuBu80}
K.\,Burrage, J.C.\,Butcher. Non-linear stability of a general class of
differential equation methods. {\em BIT} {\bf 20} (1980) 185--203.

\bibitem{Bu87}
J.C.\,Butcher. {\em The Numerical Analysis of Ordinary Differential
Equations. Runge-Kutta and General Linear Methods}, John Wiley, New York,
1987.

\bibitem{Bu93}
J.C.\,Butcher. Diagonally-implicit multi-stage integration methods. {\em
Appl. Numer. Math.} {\bf 11} (1993) 347--363.

\bibitem{Bu06}
J.C.\,Butcher. General linear methods. {\em Acta Numerica} {\bf 15} (2006)
157--256.

\bibitem{BuCaJa97}
J.C.\,Butcher, P.\,Chartier, Z.\,Jackiewicz. Nordsieck representation of
DIMSIMs. {\em Numer. Alg.} {\bf 16} (1997) 209--230.

\bibitem{BuJa93}
J.C.\,Butcher, Z.\,Jackiewicz. Diagonally implicit general linear methods
for ordinary differential equations. {\em BIT} {\bf 33} (1993) 452--472.

\bibitem{BuJa96}
J.C.\,Butcher, Z.\,Jackiewicz. Construction of diagonally implicit general
linear methods of type 1 and 2 for ordinary differential equations. {\em
Appl. Numer. Math.} {\bf 21} (1996) 385--415.

\bibitem{BuJa98}
J.C.\,Butcher, Z.\,Jackiewicz. Construction of high order diagonally
implicit multistage integration methods for ordinary differential equations.
{\em Appl. Numer. Math.} {\bf 27} (1998) 1--12.

\bibitem{BuJa04}
J.C.\,Butcher, Z.\,Jackiewicz. Construction of general linear methods with
Runge-Kutta stability properties. {\em Numer. Alg.} {\bf 36} (2004) 53--72.

\bibitem{BuJaMi97}
J.C.\,Butcher, Z.\,Jackiewicz, H.D.\,Mittelmann. Nonlinear optimization
approach to the construction of general linear methods of high order. {\em
Jour. Comput. Appl. Math.} {\bf 81} (1997) 181--196.

\bibitem{BuWr03}
J.C.\,Butcher, W.M.\,Wright. The construction of practical general linear
methods. {\em BIT} {\bf 43} (2003) 695--721.

\bibitem{IaMa} F.\,Iavernaro, F.\,Mazzia. Solving Ordinary
Differential Equations by Generalized Adams Methods: properties and
implementation techniques. {\em Appl. Num. Math.} {\bf 28}
(1998) 107--126.

\bibitem{La93} J.D.\,Lambert. {\em Numerical Methods for Ordinary
Differential Equations}, John Wiley \& Sons, Chichester, 1993.

\bibitem{MaSeTr06} F.\,Mazzia, A.\,Sestini, D.\,Trigiante.
B-spline multistep methods and their continuous extensions. {\em SIAM
Jour. Numer. Anal.}  {\bf 44} (2006) 1954--1973.

\bibitem{MaSeTr06-1} F.\,Mazzia, A.\,Sestini, D.\,Trigiante. Bs linear
multistep methods on non uniform mesh. {\em Jour. Numer.
Anal., Industrial and Appl. Math.} {\bf 1} (2006) 91--112.

\bibitem{Nord} A.\,Nordsiek. On Numerical Integration of Ordinary
Differential Equations. {\em Math. Comp.} {\bf 16} (1962) 22--49.


\bibitem{PoWeSc06}
H.\,Podhaisky, R.\,Weiner, B.A.\,Schmitt. Linearly-implicit
two-step methods and their implementation in Nordsieck form. {\em
Appl. Numer. Math.} {\bf 56} (2006) 374--387.

\bibitem{RaHi93} K.\,Radhakrishnan, A.\,Hindmarsh. Description and
Use of LSODE, the Livermore Solver for Ordinary Differential
Equations. {\em LLNL} Report UCRL-ID-113855, 1993.

\bibitem{bim} Codes {\tt BiM} and {\tt BiMD} homepage:\\
{\tt http://www.math.unifi.it/\~{}brugnano/BiM}

\bibitem{TestSet} Test Set for IVP Solvers:
{\tt http://pitagora.dm.uniba.it/\~{}testset}

\end{thebibliography}
\end{document}